\pageno=1
\input amstex
\documentstyle{amsppt}
\def\E{{\bold E}}
\def\brg{\text {br}\,\Gamma}
\magnification=1200
\topmatter
\rightline{{\it Revised 29 Jan. 1991}}
\title{
Random Walk in a Random Environment \\
and First-Passage Percolation on Trees}
\endtitle
\author
Russell Lyons\\ Robin Pemantle
\endauthor
\affil{
Department of Mathematics \\
Stanford University     \\
Stanford, CA  94305--2125     \\
Department of Statistics\\
University of California\\
Berkeley, CA  94720}
\endaffil
\thanks{The research of both authors was partially supported by NSF Postdoctoral
Fellowships. The first author was also supported by an Alfred P. Sloan
Foundation Research Fellowship.}
\endthanks
\endkeywords
\endsubjclass
\abstract{We show that the transience or recurrence of
a random walk in certain random environments on an
arbitrary infinite locally finite tree is determined
by the branching number of the tree, which is a measure
of the average number of branches per vertex.  This generalizes
and unifies previous work of the authors. It also shows that the point of
phase transition for edge-reinforced random walk is likewise determined by
the branching number of the tree. Finally, we show that
the branching number determines the rate of
first-passage percolation on trees, also known
as the first-birth problem.  Our techniques depend
on quasi-Bernoulli percolation and large deviation results.}\endabstract
\endtopmatter
\document
\def\br{\text{br}}
\def\oversigma{\raise .1truept
\hbox{\ooalign{{\raise 4.5pt\hbox{$\scriptstyle
\leftarrow$}}\crcr
{\hbox{$\textstyle\sigma$}}}}}

\def\overoversigma{\raise .2truept\hbox{\ooalign{{\raise 6.5pt\hbox
{$\scriptstyle\leftarrow$}}\crcr
{\raise 4.5pt\hbox{$\scriptstyle\leftarrow$}}\crcr
{\hbox{$\textstyle\sigma$}}}}}
\def\myarrow{\mathchoice{\displaystyle\leftarrow}{\textstyle\leftarrow}
        {\scriptstyle\leftarrow}{\scriptscriptstyle\leftarrow}}
\def\mysigma{\mathchoice{\displaystyle\sigma}{\textstyle\sigma}
        {\scriptstyle\sigma}{\scriptscriptstyle\sigma}}
\def\ovrsigma{\buildrel\myarrow\over\mysigma{}}

\overfullrule=0pt

\subheading{\S1. Introduction} 
A random walk on a tree (by which we always mean an infinite,
locally finite tree) is a Markov chain whose state space is 
the vertex set of the tree and for which the only allowable transitions 
are between neighboring vertices.  We assume throughout that all 
transition probabilities are nonzero.  For a fixed tree, the transition
probabilities may be taken as random variables, in which case
the resulting mixture of Markov chains is called ``Random Walk 
in a Random Environment (RWRE)".  The first theorem proved in this paper
is conceptually the ``least upper bound" of two previous results 
obtained by the authors (separately) about RWRE on trees.  The 
notation necessary to describe this is as follows.

Choose an arbitrary vertex as the root and let $\sigma$ be any
other vertex.  Let $\oversigma$ denote the first vertex on 
the shortest path from sigma to the root.  If $\sigma$ is at
distance at least two from the root, define $A_\sigma$ as
the transition probability from $\oversigma$ to $\sigma$ 
divided by the transition probability from $\oversigma$ to
$\overoversigma$.  We assume the following uniformity
in our random environment: all but finitely many of the random
variables $A_\sigma$ are identically distributed.  (Since the
values of $A_\sigma$ are determined by the transition probabilities in 
a way that depends on the choice of root, it may appear that
whether this condition is satisfied depends on the choice of
root, but actually a different choice of root changes only 
finitely many of the $A_\sigma$'s.)  Let $A$ denote a random
variable with this common distribution.  

By the zero-one law, a RWRE is a.s. transient or a.s. recurrent.
We shall determine the phase transition boundary; we do not
know in general when the cases on the boundary are transient or 
recurrent, but examples indicate that there is no simple rule here.
Intuitively, the larger a tree is, the more likely
a RWRE is to be transient, not only because the root
is harder to find again, but also because there will be more
branches along which the values of $A_\sigma$ are atypically
large.  This can be quantified by a large deviation calculation
which, one way or another, is behind each result in this paper.
The {\it branching number} of a tree, $\Gamma$, denoted 
$\text {br} (\Gamma)$, is a real number greater than or equal to one
that measures the average number of branches per vertex of the
tree [{\bf L1}]; the precise definition is given in \S 2. 
In [{\bf L1}, Theorem 6.6 and Remark 2],
it is shown that when $A \leq 1$, the RWRE is transient or
recurrent according to whether ${\bold E} [A] \cdot \text{br} (\Gamma )$
is greater or less than one.  In [{\bf P}, Theorem 2], it is shown
that if $\Gamma$ is a homogeneous tree or the genealogical tree
of a Galton-Watson process on (a subset of full measure of)
the event of nonextinction, then RWRE is 
a.s.~transient or a.s.~ recurrent according to whether
$p \cdot \text{br} (\Gamma )$
is greater or less than one, where $p$ is a function of the
distribution of $A$ (to be defined in the next section of this 
paper) and is equal to ${\bold E} [A]$ in the case where $A \leq 1$.  Theorem 1 
of this paper is that for any tree and any distribution of $A$,
RWRE is transient or recurrent according to whether 
$p \cdot \text{br} (\Gamma )$ is greater or less than one.  The solution
of the general case combines and simplifies techniques from
[{\bf L1}] and [{\bf P}], and provides as well a simpler expression for $p$.
As shown in [{\bf P}], edge--reinforced random walk (RRW) on an
arbitrary tree is
equal in law to a RWRE of the type discussed here.  Hence, the present
results show too that the phase transition for RRW occurs at a point
depending only on the branching number of the tree.  Our methods also
resolve the boundary case left open in [{\bf P}].

Our results on RWRE have equivalent formulations for flows in 
random electrical or capacitated networks.  In fact, that will
be important for our solution.  Such problems are often regarded as part of
percolation theory.  The second main problem we consider is also part of
percolation theory and it illustrates again how the crude behavior of
probabilistic processes on trees often involves minimal interaction between
the random variables involved and the tree structure as, moreover, the
latter enters only as a single number, {\it viz.},
the branching number of the tree.

Indeed, for our second problem, choose positive i.i.d. random variables
for each edge of a tree, $\Gamma$, regarded as transit times 
from one end to the other.  We shall give the a.s.~ rate of
fastest possible transit from any point to infinity; this is
similar to the usual problem of first-passage percolation
(or first birth) and is probably the appropriate formulation
for this setting.  Since this problem explicitly asks for
the largest deviation from mean behavior, it is not surprising
that the techniques used to find the phase transition for RWRE
also solve this problem.  Indeed, Theorem 4 gives the a.s.~
fastest transit rate to infinity as $1/m_1 (1 / \text{br} (\Gamma))$, where
$m_1$ is an inverse to a rate function $m$ defined in Section 3.
A further connection between the two problems is that 
the proofs of Theorems 1 and 4 both require results on
quasi-Bernoulli percolation [{\bf L2}] which show that a configuration
of values of $A_\sigma$ in an appropriate range, once shown
to be common enough, must percolate in an appropriate sense.
Theorem 4 also gives information concerning the asymptotic profile
of the transit times to large distances.  

We are grateful to Persi Diaconis for having introduced us to each other.

\subheading{\S2. RWRE} Given a tree, $\Gamma$, designate one of
its vertices as the root, 0.  If $\sigma$ is a vertex,
we write $|\sigma|$ for the number of edges on the shortest
path from 0 to $\sigma$.  For vertices
$\sigma$ and $\tau$, we write $\sigma\le\tau$ if $\sigma$
is on the shortest path from $0$ to $\tau$, $\sigma<\tau$
if $\sigma\le \tau$ and $\sigma\ne\tau$, and $\sigma\to\tau$
if $\sigma\le\tau$ and $|\tau|=|\sigma|+1$; in this last situation,
we call $\tau$ a {\it successor} of $\sigma$.  If $\sigma\ne 0$,
then $\overset\leftarrow\to \sigma$ denotes, as in \S1,
the vertex such that $\overset\leftarrow \to {\sigma} \rightarrow\sigma$.
The edge {\it preceding} $\sigma$, from $\overset\leftarrow\to\sigma$
to $\sigma$, is denoted $e(\sigma)$.  A {\it cutset} $\Pi$
is a finite set of vertices not including 0 such that
every infinite path from 0 intersects $\Pi$ and such that there
is no pair $\sigma,\tau\in\Pi$ with
$\sigma<\tau$.  The {\it branching number} of
$\Gamma$ [{\bf L1}] is defined by
$$\text{br}\,\Gamma:=\inf\left\{\lambda>0;\inf\limits_{\Pi}
\sum_{\sigma\in\Pi}\lambda^{-|\sigma|}=0\right\}\;.$$
The branching number is a measure of the average number of
branches per vertex of $\Gamma$.
It is  less than or equal to $\underline{\lim}_{n\to\infty} 
M_n^{1/n}$, where $M_n:=\text{card}\left\{\sigma\in\Gamma;|\sigma|=n\right\}$,
and takes more of the structure of $\Gamma$ into account
than does this growth rate.
For sufficiently regular trees, such as homogeneous trees or, more generally,
Galton-Watson trees, $\text{br}\, \Gamma = \lim_{n\to\infty} 
M_n^{1/n}$ [{\bf L1}].  

Given a random environment and the r.v.'s $A_\sigma$ as
described in \S1, we shall assume without loss of generality that 
{\it all} $A_\sigma$  are identically distributed.  Choose
further i.i.d.  $A_\sigma$  for  $|\sigma| = 1$ and set

$$C_\sigma:=\prod_{0<\tau\le \sigma} A_\sigma\;.$$
\noindent
Consider an electrical network formed from  $\Gamma$  with
conductance  $C_\sigma$  along the edge  $e(\sigma)$.  The transition
probability from  $\sigma$  to  $\tau$  is recovered by dividing the
conductance of the edge joining  $\sigma$  to  $\tau$  by the sum of
the conductances of all edges incident to  $\sigma$.  Actually, this
may not be true for  $|\sigma|\leq 1$, but we may ignore this insofar
as our interest lies in the type of the random walk.

We shall use the fact that our random walk is transient
iff the electrical network has positive conductance from $0$
to infinity (see [{\bf KSK}, Proposition 9--131]).  This, in turn, is closely
related to the question of whether the capacitated
network (with, say, water flowing instead of electricity)
formed by $\Gamma$ with (channel) capacity $C_\sigma$
through $e(\sigma)$ admits flow to infinity.  
In particular, if no water flows, then
no current flows.  (To see this, note that if the electrical
conductance is positive, then a unit potential imposed between
the root and infinity induces a current flow that is bounded
by $C_\sigma$ on $e(\sigma)$ for each $\sigma$ and is hence
an admissible water flow.)
Moreover, the converse to this is almost true, as made precise
in the proof of (1) below.  For more details, see [{\bf L1}].

\proclaim{Theorem 1} Consider a random environment on a tree
$\Gamma$  as described above with  $0 < A < \infty$ a.s.  
Let $p:=\min_{0\le x\le 1}{\bold E}[A^x]$.

\roster
\item If $p\cdot\text{br}\,\Gamma>1$, then the RWRE is
a.s.~ transient, the electrical network has positive
conductance a.s., and the capacitated network
admits flow a.s.
\item If $p\cdot\text{br}\,\Gamma <1$, then the RWRE is a.s.~ recurrent,
the electrical network has zero conductance a.s., and the
capacitated network admits no flow a.s.  More generally,
it suffices that $\inf_\Pi\sum_{\sigma\in\Pi}p^{|\sigma|}=0$.
\item If $p\cdot\overline{\lim}_{n\to\infty}M_n^{1/n}<1$,
then the RWRE is a.s.~ positive recurrent.  More generally,
it suffices that $\sum_{\sigma\in\Gamma}p^{|\sigma|}<\infty$.
\endroster
\endproclaim

\noindent{\bf Remark 1.} One might expect a part (4)
to this theorem,
as there is when $A$ is a.s.~ constant [{\bf L1}].
  However, the following example shows that 
$p\cdot\overline{\lim}
M_n^{1/n}>1$ does {\it not} imply that the RWRE is a.s.~
null recurrent or transient.  Let $\Gamma$ be a single infinite branch, to
which has been added $2^{|\sigma |+1} -1$ successors of each $\sigma$.  Each
of the added nodes has no successor.  Then $M_n^{1/n} = 2$ for all $n$. 
It is easily shown that the 
random walk will be a.s.~ positive recurrent or not according
to whether the geometric mean of $A$ is less or greater than $1/2$.  
Since we can choose $A$ to have geometric mean less than
$p$, there is no (weak)
 converse to (3).  Similar examples exist
even on trees with every vertex having at least two neighbors.

\noindent{\bf Remark 2.} The behavior on the phase transition
boundary itself will be discussed
after the proof of the theorem.

\noindent{\bf Remark 3.} Our assumptions of independence concerning
the random environment are stricter than necessary.
For example, the same proof is valid in the situation where
for some root, 0, the environments at $\sigma$ and $\tau$
are independent whenever $|\sigma|\ne |\tau|$.

\noindent{\bf Remark 4.} In case we wish to allow the
transition probabilities to be zero with
probability in $]0{,}1[$, then we must make the following changes to
the statement of the theorem.  First of all, in calculating  $p$, we
use the conventions that  $0^0 = 0$  and  $\infty^0 = 1$.  Secondly,
we change  ``a.s.'' in (1)
to ``with positive probability'';  alternatively,
instead of starting the RWRE at the root, we can say in case (1)
that there is a.s.~{\it some} vertex in $\Gamma$ at which the
RWRE is transient and from  which conductance to infinity
is positive.  To see this, suppose first that the transition probabilities,
$p_{\sigma,\tau}$  satisfy
$p_{\ovrsigma{,}\sigma} > 0$ a.s.~and ${\bold P}[p_{\sigma{,}\ovrsigma}
= 0] > 0$  for  $\sigma\not= 0$.  In this case,  
${\bold P}[A=\infty] > 0$, so that  $p = 1$  and we must show 
that the RWRE is a.s.~transient.  Indeed, it
is immediately apparent that this is the case for any subtree
consisting of a single infinite branch, hence that it is true for the
whole tree. 

More generally, now, without changing the law of the random 
environment, we may assume that the transition probabilities,
$p_{\sigma,\tau}$, have the form
$p_{\sigma,\tau}=b_{\sigma,\tau}\bar p_{\sigma,\tau}$,
where $b_{\sigma,\tau}$ takes only the values 0 and 1 
and, in fact,  $b_{\sigma,\ovrsigma}\equiv 1$; $\bar
p_{\ovrsigma{,}\sigma}$  is never  $0$;
$\{b_{\sigma,\tau}\}$ are jointly independent
of $\{\bar p_{\sigma,\tau}\}$; and $b_{\sigma,\tau}\;[\bar p_{\sigma,\tau}]$
is independent of $b_{\rho,\psi}\;[\bar p_{\rho,\psi}]$ for $\sigma\ne \rho$.
Thus, the random environment $\{p_{\sigma,\tau}\}$ can be
considered as the compound process of percolation
via $\{b_{\ovrsigma{,}\sigma}\}$ followed by use of the transition
probabilities $\{\bar p_{\sigma,\tau}\}$.  If
$\bar A$ corresponds to $\{\bar p_{\sigma,\tau}\}$
and $q:=\bold{E}[b_{\ovrsigma{,}\sigma}]$, then
$\bold{E}[A^x]=q\bold{E}[\bar A^x]$,
whence $p\cdot \text{br}\,\Gamma=\min_{0\le x\le 1}\bold{E}[\bar A^x]
(q\cdot \text{br}\,\Gamma)$.  We now combine
Theorem 1 with the fact that percolation via $\{b_{\ovrsigma{,}\sigma}\}$
leaves subtrees, the supremum of whose branching number is a.s.
$q\cdot \text{br}\,\Gamma$  [{\bf L1}, Corollary 6.3]; here, we
interpret a branching number less than $1$ to mean that the tree is finite. 
\medskip
The proof of our theorem depends on the Chernoff--Cram\'er Theorem.  We have
not found this theorem stated in the literature in the form and generality
which we require, so we state it here.  The reader may check that it
follows from the material in [{\bf D}, Chap. 1, Section 9] by a truncation
argument.  Note that we make no assumptions on the existence even of
${\E}[X]$.

\proclaim{The Chernoff--Cram\'er Theorem} Let $X$ be a real-valued random
variable and define
$$\phi(\theta):=\E[e^{\theta X}] \quad \text{and} \quad \gamma(a):=
\inf_{\theta \ge 0} (-a \theta + \log \phi (\theta)).$$
If $S_n$ denotes the sum of $n$ independent copies of $X$, then for all $a \in
{\bold R}$, the quantity
$${1 \over n}\log {\bold P}[S_n \ge na]$$ 
approaches $\gamma(a)$ from below (though not necessarily monotonically)
as $n \to \infty$.
\endproclaim

Note that in the following lemma, as well as in the proof of the theorem,
no exceptions need be made when ${\bold E} [A^x] = \infty$ for some $x$.

\proclaim{Lemma} $\min_{0 \le x \le 1} 
{\bold E} [A^x] = \max_{0 < y \le 1} \; \inf_{x \ge 0}
y^{1-x} {\bold E} [A^x] \; $.
\endproclaim

\demo{Proof} This can be verified directly by a case analysis of the point $x_0$
where ${\bold E} [A^x]$ is minimum, but it also 
follows immediately from Fenchel's
duality theorem [{\bf R}, Theorem 31.1]:  Let $f(x) := \log {\bold E}
[A^x]$ for $x \ge 0$ and $f(x) := +\infty$ for $x < 0$. 
Then $f(x)$ is convex by H\"older's inequality and lower
semicontinuous by Fatou's lemma.   Let $f^* (r)
:= \sup_x (rx - f(x))$ be the convex
conjugate of $f$.  Similarly, let $g(x)$ be the concave function that is 0 
for $x \le 1$ and $-\infty$ elsewhere, and let 
$$g^* (r)
:= \inf_x \, (rx - g(x))
= \cases
r &\text{ if }r \le 0,\cr
-\infty &\text{ if }r > 0 \cr
\endcases
$$
be its concave conjugate.  Then Fenchel's theorem asserts that $\inf_x 
(f(x) - g(x)) = \max_r (g^* (r) - f^* (r))$,
 which is the same as the statement of the lemma.  Indeed,
$\inf_x(f(x) - g(x)) = \log \min_{0\leq x\leq 1} {\bold E}[A^x]$ and
$$
\max_r(g^*(r) - f^*(r)) = \max_{r\leq 0}(r - f^*(r)) = \log \max_{0<
y\leq 1}\inf_{x\ge 0} y^{1-x}{\bold E}[A^x].\quad\qed
$$
\enddemo

\demo{Proof of Theorem 1} The assertions will be demonstrated in reverse order.

(3) We shall use the fact that the random walk is positive recurrent iff
the electrical conductances have finite sum [{\bf KSK}, Proposition 9--131]. 
Suppose that $\sum_{\sigma\in\Gamma}p^{|\sigma|}<\infty$
and that $p={\bold E}[A^x],\;0<x\le 1$.  Then
$$
{\bold E}\left[\sum_{0\ne\sigma\in\Gamma}C^x_\sigma\right]=
\sum_{\sigma\ne 0}{\bold E}\left[\prod_{0<\tau\le\sigma}A_\tau^x\right]
=\sum_{\sigma\ne 0}\prod_{0<\tau\le\sigma} {\bold E}\left[A_\tau^x\right]
=\sum_{\sigma\ne 0} p^{|\sigma|}<\infty\;,
$$
whence $\sum_{\sigma\neq 0}C_\sigma^x<\infty$ a.s. 
In particular, $C_\sigma< 1$ for all but finitely many $\sigma$
a.s. Since $C^x_\sigma\ge C_\sigma$ for $C_\sigma<1$,
it follows that $\sum_{\sigma\neq 0}C_\sigma<\infty$
a.s.  

(2) Here we use the fact that the random walk is recurrent iff no 
electrical current flows.  
Suppose that $\inf_\Pi\sum_{\sigma\in\Pi}p^{|\sigma|}=0$
and that $p= {\bold E}[A^x],\;0<x\le 1$.  Then as above,
$$
{\bold E}\left[\sum_{\sigma\in\Pi}C_\sigma^x\right]=
\sum_{\sigma\in\Pi}p^{|\sigma|}\;,
$$
whence if $\sum_{\sigma\in\Pi_n}p^{|\sigma|}\to 0$,
we have, as above,
$$\underline{\lim}_{n\to\infty}\sum_{\sigma\in\Pi_n}C_\sigma=0
\;\text{a.s.}$$
by virtue of Fatou's lemma.  By (the trivial half of) the max-flow min-cut
 theorem, the capacitated network admits no
 flow a.s.
  Hence no electrical current flows, and the random
walk is a.s.~recurrent. 

(1) This part uses the fact that if water flows even when $C_\sigma$ is
reduced exponentially in $|\sigma |$, then electrical current flows
and the random walk is transient [{\bf L1}, Corollary 4.2].    
If $p \cdot \text{br}\Gamma>1$, let $y\in ]0,1]$ be such that
$p=\inf_{x\in\bold R} y^{1-x}{\bold E}[A^x]$.
By the Chernoff--Cram\'er Theorem, there exists $k\ge 1$
such that for $|\sigma|=k$,
$${\bold P}\left[C_\sigma\ge y^k\right]>(y\cdot\text{br}\Gamma)^{-k}\;.$$
Let $\epsilon>0$ be sufficiently small that for $|\sigma|=k$,
$$q:={\bold P}[C_\sigma\ge y^k\;\&\;\forall\; 0<\tau\le \sigma\;
A_\tau\ge\epsilon]>(y\cdot\text{br}\,\Gamma)^{-k}\;.$$
Let $\Gamma^k$ be the tree whose vertices are
$\left\{\sigma\in\Gamma;k\Bigm||\sigma|\right\}$ such that
$\sigma\to\tau$ in $\Gamma^k$ iff $\sigma\le\tau$
and $|\sigma|+k=|\tau|$ in $\Gamma$.
It is easily verified that  br$\,\Gamma^k=(\text{br}\,\Gamma)^k$.

Form a random subgraph, $\Gamma^k(\omega)$, of $\Gamma^k$
by deleting those edges $\sigma\to \tau$ where
$$\prod_{\scriptstyle \sigma<\rho\le\tau\atop\scriptstyle
\rho\in\Gamma} A_\rho<y^k\;\text{or}\;\exists\rho\in\Gamma\;
(\sigma<\rho\le\tau\;\&\;A_\rho<\epsilon)\;.$$
This is an edge percolation on $\Gamma^k$
for which each edge is present with probability $q$ 
 and the presence of  edges with distinct ``preceding''
vertices are mutually 
independent.
  In particular, it is a quasi-Bernoulli percolation process on
$\Gamma^k$ (see [{\bf L2}] for the general definition of quasi-Bernoulli 
percolation). 
Choose $w\in ](yq^{1/k}\text{br}\,\Gamma)^{-1},1[$.
 Since $q\cdot\text{br}\,\Gamma^k>(wy)^{-k}> 1$,   
there is almost surely
a subtree $\Gamma^*$ of $\Gamma^k (\omega )$, not necessarily beginning
at the root, that has branching number larger than $(wy)^{-k}$
(combine the method of proof of Corollary 6.3 or Proposition 6.4 of
[{\bf L1}] with Theorem 3.1 of [{\bf L2}]).  
Any subtree $\Gamma^*$ of $\Gamma^k$ induces a subtree $\Gamma^\prime$ 
of $\Gamma$ 
whose vertices are those $\rho\in \Gamma$ for which $\exists\sigma,\tau\in
\Gamma^*$ such that $\sigma\le \rho\le \tau$.
Thus, there is almost surely a subtree
$\Gamma^\prime$ of $\Gamma$ induced by $\Gamma^k (\omega )$ with the 
following three properties: $\text{br}\, \Gamma^\prime >(wy)^{-1}$;
 $A_\sigma \ge \epsilon$
along each edge; and
for every $\sigma,\tau\in\Gamma^\prime$ such that $k
\Bigm| |\sigma|=|\tau|-k$, we have $\prod_{\sigma<\rho\le\tau}
A_\rho\ge y^k$.  It follows that 
$$\inf_{\Pi^\prime}\sum_{\sigma\in\Pi^\prime}w^{|\sigma|}C_\sigma
\ge\inf_{\Pi^\prime}\epsilon^{k-1}\sum_{\sigma\in\Pi^\prime}(wy)^{|\sigma|}
>0\;,$$
where $\Pi^\prime$ is any cutset of $\Gamma^\prime$.
By [{\bf L1}, Corollary 4.2], $\Gamma^\prime$ has positive
conductance, hence so does $\Gamma$ (a.s.).
We may now deduce that the RWRE is transient a.s.~and
that the capacitated network admits flow 
a.s.~(for example, the current flow
from an imposed
unit potential from 0 to infinity). \qed
\enddemo

\noindent{\bf Remark 5.} The proofs of (1) and (2) use the different
expressions of $p$ given by the lemma.  Since the equality of the two
expressions is not intuitive, the fact that  (1) and (2) meet (cover all
possibilities except $p \cdot \text{br} \, \Gamma = 1$) may seem miraculous.
The following discussion is intended to explain ``why" (1) and (2) meet.  Part (1) is true because there
are enough branches of the tree along which the geometric mean of the 
$A_\sigma$'s
exceeds a certain value, $y$, in order to force transience.  This is shown 
directly from the second expression for $p$, which is just a large deviation
rate.  For each $y$ individually, there is a converse to this, which is that
the conductances to those parts of a cutset where the geometric mean of 
the $A_\sigma$'s
is close to $y$ approach zero as the cutset gets far from the root.
This uses the other half ({\it i.e.}, the upper bound) of the
Chernoff--Cram\'er Theorem not used in (1).
It is possible to show that the total conductance to
a cutset goes to zero by ``integrating'' this fact 
over $y$.  Such an approach
can be used to derive (2) from the second expression for $p$, but it involves
several calculations (see [{\bf P}] where this is done for positive 
recurrence).  These are avoided by using
 the simpler expression for $p$ provided by the lemma.
\medskip
We now turn our attention to the behavior of the random walk on the
phase transition boundary,  $p\cdot\text{br}\,\Gamma = 1$.
Here, the type of the walk depends on further structure of the tree and
of $A$.
 First,
we remark that even when  $A$  is constant, the walk may be either
transient or recurrent [{\bf L1}].
\noindent
 Furthermore, there are trees, $\Gamma$, for which $p =\text{br}
(\Gamma) = 1$
 and RWRE is a.s.~transient, yet simple random walk on $\Gamma$
 ({\it i.e.}, $A = p = 1$ almost surely) is recurrent. On the other hand, 
 there are trees, $\Gamma$, for which $\br (\Gamma) > 1$ and for which 
 there are a.s.~recurrent RWRE's with $p = \br(\Gamma)^{-1}$ , even though
 the RWRE with the deterministic environment $A = p$ is transient.
We do not know whether the difference in behavior of the walk for random
as compared to deterministic environments always depends, as above, on $\br\,
\Gamma$.
 We hope to clarify this behavior at a later date, but for now,
 suffice it to record the following extension to Part (2) of Theorem 1,
 which covers many of the boundary cases, $p \cdot \br(\Gamma) = 1$.

\proclaim{Proposition 2}  Under the same hypotheses as in Theorem 1,
if  $\underline{\lim}_{\Pi\to\infty}\sum_{\sigma\in \Pi}p^{|\sigma|} < \infty$
(in other words, if there are cutsets  $\Pi_n$  such that
$\inf\{|\sigma|; \sigma\in\Pi_n\}\to\infty$  and
$\sup_n\sum_{\sigma\in\Pi_n} p^{|\sigma |} < \infty $), then the RWRE
is a.s.~recurrent and the electrical network has zero conductance a.s.
\endproclaim

\demo{Proof}  Let  $\Pi_n$  be as indicated.  We have
$$
{\bold E}[\sup_n \sum_{\sigma\in\Pi_n} C^x_\sigma]\leq \sup
{\bold E}[\sum_{\sigma\in\Pi_n} C^x_\sigma] = \sup \sum_{\sigma\in\Pi_n}
p^{|\sigma|} < \infty,
$$
\noindent whence
$$
\sup \sum_{\sigma\in \Pi_n} C_\sigma \leq
\sup(\sum_{\sigma\in\Pi_n}C^x_\sigma)^{1/x} < \infty \;\text{a.s.}
$$
This implies a.s.~recurrence by virtue of [{\bf L1}, Corollary
4.2].   \qed 

The capacitated network may admit flow in these circumstances, as the
case where  $\Gamma$  is a binary tree and  $A\equiv 1/2$  shows.  A
complete answer for both electrical and capacitated networks may be
given when  $\Gamma$  is homogeneous or produced by a Galton-Watson
process.  The following theorem simplifies, refines and extends
Theorem 2 of [{\bf P}].  (The assumptions there that the progeny
distribution be bounded and that ${\bold E}[\log A]$  exist are now
seen to be unnecessary.)
\enddemo

\proclaim{Theorem 3}  Let  $\Gamma$  be the genealogical tree of a
Galton-Watson branching process with mean  $m > 1$.  Consider a
random environment as in Theorem 1.

\roster
\item If  $pm > 1$, then given nonextinction, the RWRE is a.s.~
transient, the electrical network has positive conductance a.s., and
the capacitated network admits flow a.s.
\item If  $pm \leq 1$, then given nonextinction, the RWRE is
a.s.~recurrent, the electrical network has zero conductance a.s.,
and, unless both  $A$ and  $\Gamma$  are constant, the capacitated
network admits no flow a.s.
\item If  $pm < 1$, then the RWRE is a.s.~positive recurrent.
\endroster
\endproclaim

\noindent {\bf Remark 6.} When  $pm = 1$, the RWRE may be either null or
positive recurrent  [{\bf L3}, Theorem 3.2].

\demo{Proof} Part (1) follows immediately 
from Theorem 1, since given nonextinction, the
genealogical tree has branching number $m$  a.s. [{\bf L1},
Proposition 6.4].  The first two parts of (2) follow from the same
argument used in the proof of Proposition 2:

$$
{\bold E}[\sum_{|\sigma|= n} C^x_\sigma] = m^n p^n\leq 1,
$$
\noindent whence

$$
\sup_n \sum_{|\sigma|= n} C_\sigma < \infty\;\text{a.s.}
$$
\noindent Likewise, the proof of (3) is similar to the proof of Part
(3) of Theorem 1.  It remains to establish the assertion in (2) on
the capacitated network when  $pm = 1$.

For $0 < t \leq 1$, let $F^{(t)}$ be the maximum flow to
infinity in the network on $\Gamma$ with capacities $C_\sigma^t$
along $e(\sigma)$.  Thus, we are interested in showing that
$F^{(1)} = 0$ almost surely.  Now as in Proposition 2 above,
$$ (F^{(1)} )^x = \inf_\Pi \left ( \sum_{\sigma \in \Pi} C_\sigma 
     \right )^x \leq \inf_\Pi \left ( \sum_{\sigma \in \Pi} C_\sigma^x
     \right ) = F^{(x)} , $$
so that it suffices to show that $F^{(x)} = 0$ almost surely.
For $|\sigma| = 1$, let $F_\sigma^{(x)}$ be the maximum flow in
the subtree $\{ \tau \in \Gamma \,;\, \sigma \leq \tau \}$ with capacities
$C_\tau^x / A_\sigma^x$.  Thus, $F_\sigma^{(x)}$ has the same law as
$F^{(x)}$ does.  It is easily seen that 
$$F^{(x)} = \sum_{|\sigma| = 1} A_\sigma^x (1 \wedge F_\sigma^{(x)} ). $$
Taking expectations yields
$$
\E [ F^{(x)} ] = mp \E [ 1 \wedge F^{(x)} ] = \E [ 1 \wedge F^{(x)} ] .
$$
Therefore $F^{(x)} \leq 1$ almost surely.  In addition, we have, by
independence,
$$\parallel F^{(x)} \parallel_\infty \;=\; \parallel \sum_{|\sigma| = 1}
    A_\sigma^x \parallel_\infty \cdot \parallel F^{(x)} 
    \parallel_\infty .$$
Since $\parallel \sum_{|\sigma| = 1} A_\sigma^x \parallel_\infty > 1$
unless both $A$ and $\Gamma$ are constant, this shows that
$F^{(x)} = 0$ almost surely unless both $A$ and $\Gamma$ are constant.  \qed
\enddemo

\subheading{\S3.  First-passage percolation} Given a tree
$\Gamma$ rooted at 0, the {\it boundary},
$\partial\Gamma$, of $\Gamma$ is the space of infinite
paths beginning at $0$ which go through no vertex more than
once.  This is a compact space with metric $d(s,t)=e^{-n}$, where $n$ is the 
number
of edges common to $s$ and $t$.  Changing the root gives  essentially
the same
boundary with an equivalent metric.
Suppose we are given real-valued i.i.d.r.v.'s
$X_\sigma$ for each edge $e(\sigma)$.
Let $S_\sigma=\sum_{0<\tau\le\sigma}X_\tau$.
As explained in \S1, the random variable
$$
\inf_{s\in\partial\Gamma}\overline{\lim}_{\sigma\in s}
\frac{S_\sigma}{|\sigma|}
$$
may be thought of as the reciprocal of the fastest sustainable transit rate
to infinity if  $X_\sigma > 0$.  The calculation of this rate depends on the
non-decreasing function

$$
m(y)=\inf_{x\le 0}{\bold E}[e^{x(X-y)}]\;,$$
where $X$ has the same law as every $X_\sigma$.
As the infimum of linear functions, $\log m$ is concave.  Write
$$
m_1(z) := \sup \{ y \, ; \, m(y) < z \} \, ,
$$
so that $m_1$ is a sort of inverse function to $m$.  Note that $m$
cannot be constant unless $m \equiv 1$ nor have range $\{0, \, 1\}$.  Hence $m$
is strictly increasing where ]0, 1[ -- valued by log--concavity and
$$
m_1(z) = \inf \{ y \, ; \, m(y) > z \}
$$
for $0 \le z \le 1$ except if $m \equiv 1$ and $z = 1$.

\proclaim{Theorem 4} Unless $\text{br}\,\Gamma = 1$ and $m \equiv 1$, we have
$$
\inf_{s\in\partial\Gamma}\overline{\lim}_{\sigma\in s}
\frac{S_\sigma}{|\sigma|}
=\inf\{\lim_{\sigma\in s}\frac{S_\sigma}{|\sigma|}\, ; \, s \in
\partial\Gamma\; \&\; \lim_{\sigma\in
s}\frac{S_\sigma}{|\sigma|}\;\text{exists}\}
=m_1((\text{br}\,\Gamma)^{-1}) \;\text{a.s.}$$ 
and 
$$
\dim\left\{s\in\partial\Gamma\,;\,\overline{\lim}_{\sigma\in s}
\frac{S_\sigma}{|\sigma|}\le y\right\}
=\dim\left\{s\in\partial\Gamma\,;\,\lim_{\sigma\in s}
\frac{S_\sigma}{|\sigma|} = y\right\}
=\log (m(y)\text{br}\,\Gamma)\;\text{a.s.}
$$
for $m_1((\text{br}\,\Gamma)^{-1})\le y < \sup_n {\bold E}[X\wedge n]$.
\endproclaim

\noindent{\bf Remark 7.} The statement concerning Hausdorff
dimension ({\it cf.} [{\bf L1}, \S7])
may be interpreted as giving information on the asymptotic
profile of transit times.  We are grateful to Yuval Peres for the
statements and proofs concerning ``lim" (rather than
``$\overline{\lim}$").  Peres also has examples showing that when
$\text{br}\, \Gamma = 1$ and $m \equiv 1$, the rate depends on more information.

\demo{Proof}Suppose first that
$m(y)\text{br}\,\Gamma<1$.
Then by the Chernoff--Cram\'er Theorem,
$$\inf_\Pi{\bold E}\left[\sum_{\sigma\in \Pi} \text{\bf 1}_{S_\sigma\le
y|\sigma|}\right]\le
\inf_\Pi\sum_{\sigma\in\Pi}m(y)^{|\sigma|}=0\;,$$
whence there are cutsets $\Pi_n\to\infty$
({\it i.e.}, $\min\{|\sigma|;\sigma\in\Pi_n\}\to\infty$) such that
$$\underline{\lim}_{n\to\infty}\text{card}\{\sigma\in
\Pi_n;S_\sigma\le y|\sigma|\}=0\;\text{a.s.}$$
In other words, a.s.~for infinitely many $n$,
$$\forall \sigma\in\Pi_n\; S_\sigma>y|\sigma|\;.$$
Therefore a.s.
$$
\forall s\in\partial\Gamma\;\overline{\lim}_{\sigma\in s}
\frac{S_\sigma}{|\sigma|} \ge y\;.
$$

On the other hand, if
$m(y)\text{br}\,\Gamma>1$, then by the Chernoff--Cram\'er
Theorem, there is a $k\ge 1$ such that for $|\sigma|=k$,
$$
{\bold P}[S_\sigma\le ky]>(\text{br}\,\Gamma)^{-k}\;.
$$
Let $M$ be sufficiently large that for $|\sigma|=k$,
$$
q:={\bold P}\left[S_\sigma\le ky\;\&\;\forall\; 0<\tau\le\sigma\;
X_\tau\le M\right]>(\text{br}\,\Gamma)^{-k}\;.\leqno(3.1)
$$
\noindent Let  $\Gamma^k$  be as in the proof of Theorem 1.
Form a random subgraph, $\Gamma^k(\omega)$,
of $\Gamma^k$ by deleting those edges $\sigma\to\tau$
where
$$\sum_{\scriptstyle \sigma<\rho\le\tau\atop\scriptstyle
\rho\in\Gamma} X_\rho>ky\;\text{or}\;\exists \rho\in \Gamma
\;(\sigma<\rho\le\tau\;\&\;X_\rho>M)\;.$$
This is a quasi-Bernoulli percolation process on
$\Gamma^k$ [{\bf L2}].
Since $q\text{br}\,\Gamma^k>1$, percolation occurs a.s.
That is, there is a.s.~an $s\in \partial\Gamma$ whose image
in $\partial\Gamma^k$ is, except for a finite set, contained in
$\Gamma^k(\omega)$; for such $s$,
we have $\overline{\lim}_{\sigma\in s} S_\sigma/|\sigma|\le y$.
This establishes that
$ \inf_{s\in\partial\Gamma}\overline{\lim}_{\sigma\in s}
S_\sigma/|\sigma|
=m_1((\text{br}\,\Gamma)^{-1}) \;\text{a.s.}$

Additional information can be extracted from the argument of the last
paragraph.  Embed $\Gamma$ in the upper half--plane with its root at the
origin and order $\partial\Gamma$ clockwise.  Given $y$ such that
$m(y)\,\brg > 1$, let $\Gamma^k_y$ be the tree denoted $\Gamma^k$ above
and let $s(y)$ be the minimal element of $\partial\Gamma$ whose tail lies
in $\Gamma^k_y(\omega)$.  Thus, $s(y)$ is defined a.s.  For any
$|\sigma| = k$, set
$$
\psi(y) := \E[\frac{S_\sigma}{k} \mid S_\sigma \le ky \; \& \; \forall\, 0
< \tau \le \sigma \quad X_\tau \le M].
$$
Recall that $k$ and $M$ depend on $y$.  By the strong law of large
numbers, we have
$$
\lim_{\sigma \in s(y)} \frac{S_\sigma}{|\sigma|} = \psi(y) \quad\text {a.s.}
$$
Since $\psi(y) \le y$ and $y$ is arbitrary subject only to $m(y)\,\brg >
1$, this establishes the remainder of the first assertion of the
theorem.  We claim, moreover, that for any $y < \sup_n \E[X\wedge
n] $ such that $m(y)\,\brg > 1$,
$$
\text {a.s.}\;\exists s \in \partial\Gamma \quad \lim_{\sigma \in s}
\frac{S_\sigma}{|\sigma|} = y \; .
\leqno(3.2)
$$
Indeed, given such $y$, find $k$ and $M$ such that (3.1) holds and $y <
\E[X\wedge M]$.   We may write
$$
y=\alpha \psi(y) + (1 - \alpha) \E [X \mid X \le M]
$$
for some $\alpha \in [0,\,1]$.  Choose a sequence ${\Cal N}$ of density
$\alpha$ in {\bf N} and percolate as before on $\Gamma^k$, {\it except}
that the edges preceding vertices $\sigma \in \Gamma^k$ for which $|\sigma|
\not \in {\Cal N}$ survive a.s., rather than with probability $q$.  We now
find by similar reasoning to the above that
$$
\lim_{\sigma \in s(y)}\frac{S_\sigma}{|\sigma|} = y \quad \text{a.s.},
$$
thereby validating (3.2).

The remainder of the theorem follows from general considerations.
Namely, denote the sets in question by
$$
\leqalignno{
E(y) &:= \{ s \in \Gamma \,;\, \overline{\lim}_{\sigma \in
s}\frac{S_\sigma}{|\sigma|} \le y \},\cr
F(y) &:= \{ s \in \Gamma \,;\, {\lim}_{\sigma \in
s}\frac{S_\sigma}{|\sigma|} = y \}.
\cr}
$$
By the 0--1 law, dim $E(y)$ and dim $F(y)$ are constant a.s.  As $E(y)$
and $F(y)$ are clearly Borel and $F(y) \subseteq E(y)$, [{\bf L1}, \S
7] implies that it suffices to show that if independent Bernoulli
percolation with survival parameter $p$ is performed on $\Gamma$, then
for $pm(y)\,\brg < 1$, a.s.~no point of $E(y)$ survives, while for
$pm(y)\,\brg > 1$, with positive probability some point of $F(y)$ does
survive.  These conditions in fact follow from [{\bf L1}, Corollary
6.3], what was shown above, and Fubini's theorem. \qed

\enddemo

\bigskip
\centerline{\smc References}
\medskip

\item{[{\bf D}] }Richard Durrett, {\it Probability: Theory and
Examples.} Wadsworth \& Brooks/Cole, Belmont, CA, 1991.

\item{[{\bf KSK}] }John G. Kemeny, J. Laurie Snell 
and Anthony W. Knapp, {\it Denumerable Markov Chains.} 2nd ed.
Springer, New York, 1976.

\item{[{\bf L1}] }Russell Lyons, {\it Random walks and percolation
on trees,}  Ann. Probab. {\bf 18} (1990), 931--958.

\item{[{\bf L2}] } \vrule width 1.5cm height 0.4pt depth 0.4pt ,
{\it The Ising model and percolation on trees and tree-like graphs,}
Commun. Math. Phys. {\bf 125} (1989), 337--353. 

\item{[{\bf L3}] }\vrule width 1.5cm height 0.4pt depth 0.4pt , 
{\it Random walks, capacity, and percolation on trees,} preprint, 1990.

\item{[{\bf P}] }Robin Pemantle, {\it Phase transition in reinforced
random walk and RWRE on trees,} Ann. Probab. {\bf 16} (1988),
1229--1241.

\item{[{\bf R}] }R. Tyrell Rockafellar, {\it Convex Analysis.}
Princeton University Press, Princeton, 1970.

\bigskip
\noindent {\it Current Addresses:}  

Department of Mathematics, Indiana University, Bloomington, IN  47405

Department of Mathematics, Oregon State University, Corvallis, OR  97331-4605

\enddocument